\documentclass[12pt, fullpage]{article}
\usepackage{fullpage}
\usepackage{latexsym}
\usepackage{amsthm}
\usepackage{amssymb}
\usepackage{amsmath}
\newtheorem{theorem}{Theorem}[section]

\newtheorem{lemma}[theorem]{Lemma}
\newtheorem{definition}[theorem]{Definition}
\newtheorem{conjecture}[theorem]{Conjecture}

\newcommand{\AD}{\mathcal{A}}
\newcommand{\C}{\mathbb{C}}
\newcommand{\F}{\mathbb{F}}
\newcommand{\II}{{\rm II}}
\newcommand{\bi}{{\bf i}}
\newcommand{\bj}{{\bf j}}
\newcommand{\bk}{{\bf k}}

\renewcommand{\L}{\mathcal{L}}

\renewcommand{\P}{\mathbb{P}}
\newcommand{\Q}{\mathbb{Q}}
\newcommand{\R}{\mathbb{R}}
\newcommand{\X}{\mathcal{X}}
\newcommand{\Z}{\mathbb{Z}}

\DeclareMathOperator{\Gal}{Gal}

\DeclareMathOperator{\PGL}{PGL}

\title{A genus $2$ family with $226$ sections}
\author{Genya Zaytman}
\date{}

\begin{document}

\maketitle

\section{Introduction}

Faltings' theorem \cite{Fa1},\cite{Fa2} (formerly the Mordell conjecture \cite{mor}) states that a curve of genus greater than one over any number field has only finitely many points.  Again a natural question is how many points can such a curve have.  Caporaso, Harris, and Mazur \cite{chm} have shown that the weak Bombieri-Lang conjecture implies that for any number field $F$ and any integer $g\ge2$ there is an absolute upper bound $B(F,g)$ on the number of points on a genus $g$ curve over $F$. Furthermore, the strong Bombieri-Lang conjecture implies that for each genus $g\ge2$, there is an absolute bound $C(g)$ depending on the genus --- but not on the field --- such that over any number field, only finitely many curves of genus $g$ have more than $C(g)$ points.  Again we can ask what those two bounds are and, as it turns out, it helps to consider families that come from K3 surfaces.

Specifically, we will consider the case $g=2$.  We use a K3 surface $\X$ that is a double cover of $\P^2$ ramified over a smooth sextic curve $C$, so every pencil of lines gives us a family of genus $2$ curves.  Any line that is tangent to $C$ at $3$ points will lift to a pair of curves on $\X$ that become sections of the family given by any pencil of lines.  In section \ref{genus 2} we will construct an K3 surface (over a number field) and corresponding sextic with $64$ such tritangents.  Furthermore, there will turn out to be other rational curves in $\P^2$ of higher degree that also meet $C$ only at tangent points.  By suitably choosing the pencil and performing suitable base changes we find a family of genus $2$ curves with $226$ sections.  This is the current best record, the previous \cite{elk g2} being $150$ sections, which remains the record for a family over~$\Q$.
\\
\\
\emph{Acknowledgments.}  I'd like to thank Noam Elkies for the great deal of mathematical assistance 
he has provided on this project and throughout my time in graduate school.  This research was partially supported by an National Science Foundation Graduate Fellowship.

\section{Background}

We will need to review the theory of K3 surfaces.  However, before we can do that, we must first deal with the general theory of lattices.  A much more thorough discussion of Lattices to which we will occasional refer to can be found in \cite{cs}.

\subsection{Lattices}

\begin{definition}
For our purposes a \emph{lattice} is a free module $L$ of finite rank over $\Z$ together with a symmetric nondegenerate bilinear pairing $L\times L\to\Q$, $(u,v)\mapsto u\cdot v$.  The \emph{norm} of an element $C\in L$ is defined as $N(u)=u\cdot u$.  A lattice is said to be \emph{integral} if the image of the dot product lies in $\Z$ and \emph{even} if the image of the norm map lies in $2\Z$.
\end{definition}

Notice that every even lattice is integral since by bilinearity
$$u\cdot v=\frac{1}{2}(N(u+v)-N(u)-N(v)).$$

\begin{definition}
A lattice is said to be \emph{positive definite} (resp. \emph{negative definite}) if the image of the norm map lies in $\Q^{\ge0}$ (resp. $\Q^{\le0})$ and \emph{indefinite} if it is neither positive nor negative semi-definite.
\end{definition}

We define a sublattice and in the obvious manner.

\begin{definition}
A sublattice $L'\subset L$ is said to be \emph{primitive} if the group $L/L'$ is torsion-free.
\end{definition}

\begin{definition}
The \emph{orthogonal complement} of a sublattice $L'\subset L$ is the lattice $L''=\{v\in L \;|\; v\cdot w=0\; {\rm for\; all}\; w\in L' \}$.
\end{definition}

\begin{definition}
The \emph{determinant} of a lattice is the determinant of the matrix $\{v_i\cdot v_j\}_{ij}$ where the $v_i$'s are a basis of $L$.  An integral lattice is said to be \emph{unimodular} if its determinant is $\pm1$.
\end{definition}

Note that the determinant is independent of the basis.  We will on occasion extend the pairing to a map $(L\otimes R)\times (L\otimes R)\to(\Q\otimes R)$ for other rings $R$.

An important operation on lattices is taking the dual.

\begin{definition}
The \emph{dual} of a lattice $L$, denoted $L^*$, is defined to be
$$\{u\in L\otimes\Q \;|\; u\cdot v\in\Z\; {\rm for\; all}\; v\in L \}.$$
\end{definition}

In particular if $L\subset L^*$ if and only if $L$ is integral and $L=L^*$ if and only if $L$ is unimodular.

\begin{definition}
The \emph{determinant group} of an even lattice $L$ is the group $d(L)=L^*/L$.  The \emph{determinant form} of $L$ is the map $N_d:d(L)\to(\Q/2\Z)$ defined by $N_d(C)\equiv N(C) \pmod{2\Z}$.
\end{definition}

Note that $N_d$ is well-defined because the lattice is even.  The order of the determinant group is the absolute value of the determinant.  If $L\subset L'$ is primitively embedded and $T$ is the orthogonal complement then $d(L)$ and $d(T)$ are isomorphic groups with anti-isomorphic forms \cite[Chap.~4]{cs}.

\subsection{K3 Surfaces and the Torelli Theorem}

\begin{definition}
An \emph{algebraic K3 surface} is a simply-connected algebraic surface with trivial canonical bundle.
\end{definition}

All K3 surfaces $\X$ over $\C$ have $H^2(\X,\Z)\cong\II_{3,19}$, where $H^2(\X,\Z)$ is the second topological cohomology group with cup product and $\II_{3,19}$ is the unique even unimodular lattice of signature $(3,19)$.  Furthermore, the Hodge numbers corresponding to the decomposition of $H^2(\X,\Z)\otimes\C$ into $H^{2,0}$, $H^{1,1}$, and $H^{0,2}$ are $1,\,20,\,1$.

There is a Torelli theorem for K3 surfaces stating that a K3 surface $\X$ over $\C$ is determined up to complex conjugation by the cup product and Hodge decomposition of $H^2(\X,\Z)\otimes\C$.  More precisely: let $\pi$ be the projection map of $H^2(\X,\Z)\otimes\C$ onto $H^{2,0}$.  By applying  complex conjugation we see that the kernel of $\pi$ restricted to $H^2(\X,\Z)\otimes\R$ lies entirely in in $H^{1,1}$.  Hence by dimension count $\pi|_{H^2(\X,\Z)\otimes\R}$ is surjective onto $H^{2,0}$ and $\ker \pi|_{H^2(\X,\Z)\otimes\R}\otimes\C=H^{1,1}$; furthermore, the induced  cokernel norm on $H^{2,0}$ is a positive scalar multiple of the standard norm on a $1$-dimensional complex vector space.  By a theorem of Lefschetz \cite{lef} the kernel of $\pi$ restricted to $H^2(\X,\Z)$ is the lattice of algebraic cycles, also called the N\'eron-Severi lattice, $NS(\X)$.

\begin{definition}
The sublattice of $H^2(\X,\Z)$ orthogonal to $NS(\X)$ is called the \emph{transcendental lattice} and will be denoted $T(\X)$.
\end{definition}

\begin{theorem}\label{torelli}\emph{(Torelli Theorem for K3 surfaces.)~\cite{ps},\cite{fried}}
Let $\X$ and $\X'$ be K3 surfaces.  Any isomorphism $H^2(\X,\Z)\to H^2(\X',\Z)$ that takes effective cycles to effective cycles, or equivalently that takes an ample divisor to an ample divisor, and takes the form $\pi_\X$ to a form proportional to $\pi_{\X'}$ comes from a unique isomorphism of surfaces $\X\to \X'$.
\end{theorem}

The Torelli theorem also has a converse.

\begin{theorem}\label{torelli converse}\emph{\cite{fried}}
Let $\pi:\II_{3,19}\to\C$ be a linear map whose kernel has a vector of negative norm and such that the induced map $\II_{3,19}\otimes\R\to\C$ with the standard pairings is a projection up to positive scalars.  Then there exists a K3 surface $\X$ and a map $f:\II_{3,19}\to H^2(\X,\Z)$ that takes $\pi$ to a form proportional to $\pi_\X$.  Furthermore, if $v\in \ker\pi$ has positive self intersection then we can choose $f$ so that $f(v)$ is quasi-ample.
\end{theorem}

Given a hyperbolic lattice $L$ primitively embedded in $H^2(\X,Z)$, the Torelli theorem can be used to determine the moduli space of K3 surfaces $\X$ polarized by a primitive embedding $L\to NS(\X)$ as a nonempty union of complex manifolds of dimension $20-\operatorname{rank}(L)$.  Unfortunately, the Torelli theorem is not constructive enough to write down equations for the moduli space, let alone for an explicit K3 surface.  Also since this correspondence is only defined over $\C$, it provides no arithmetic information such as the field of definition or reduction modulo primes.

\section{Constructing the genus $2$ family}
\label{genus 2}

\subsection{Families of Curves}
\label{family}

Given a K3 surface (or indeed any surface) $\X$, expressing $\X$, possibly with some points blown up, as a family of genus $g$ curves over a rational curve corresponds to finding a smooth positive divisor $D$ of genus $g$  on the surface and choosing a $2$-dimensional subspace of $\L(D)$.  Sections correspond to irreducible positive divisors $S$ that have intersection number $S\cdot D=1$.    On a K3 surface, since the canonical bundle is trivial, the adjunction formula gives $D\cdot D=\chi(D)=2g-2$, and $S\cdot S=\chi(S)=-2$.  Conversely, by Hirzebruch-Riemann-Roch if $D$ is a divisor on a K3 surface satisfying $D\cdot D=2g-2$ then either $D$ or $-D$ is a positive divisor whose line bundle has $h_0\ge g+1$.  Furthermore, if the $D$ above is positive and $S$ is a divisor satisfying $S\cdot D=1$ and $S\cdot S=-2$ then $S$ is also positive.  Unfortunately it is possible for two such $S$'s to correspond to the same section if one of them contains a component $C$ of a reducible fiber, thus reducing the actual number of sections.  Such a $C$ satisfies $C\cdot D=0$ and $C\cdot C=-2$.  If $g\ge2$, one can take the orthogonal projection $\psi$ of the N\'eron-Severi lattice from $D$ to get a negative definite lattice with $\psi(S)\cdot\psi(S)=-2-\frac{1}{2g-2}$ and $\psi(C)\cdot\psi(C)=-2$, hence one wants a lattice of the correct lattice genus (so it can be glued to $D$ to get an integral even lattice \cite[Chap.~4]{cs}) of small enough discriminant to have many vectors of norm $-2-\frac{1}{2g-2}$ but not small enough to have too many vectors of norm $-2$.   Notice also that if there are no reducible fibers and $g>1$ then $D$ has positive self-intersection and positive intersection with all positive divisors; thus by the Nakai-Moishezon criterion $D$ is an ample divisor.

Specializing to the case when $g=2$, we see that $D\cdot D=2$ and by Hirzebruch-Riemann-Roch $\ell(D)\ge 3$.  Assuming we are in the case where $H^1(\X,O(D))=0$ (or taking a generic $3$-dimensional subspace of $\L(D)$) gives us a map $\phi:\X\to\P^2$.  Since $D\cdot D=2$, this map is a double cover.  To compute the degree of the ramification divisor, consider a generic line $E\subset\P^2$.  By construction $\phi^{-1}(E)$ has genus $2$, so by Riemann-Hurwitz $\phi$ restricted to $\phi^{-1}(E)$ is ramified at $6$ points.  In other words, $E$ intersects the ramification divisor at $6$ points; therefore, the ramification divisor is a curve $C_F:F(X,Y,Z)=0$ of degree $6$.

If $C$ is a component of a fiber, i.e. $C\cdot D=0$, then $\phi(C)$ must be a point (since otherwise $C$ would intersect $D=\phi^{-1}(E)$) which must thus be a singular point of $F$.  If $S$ is a section and therefore has $S\cdot D=1$ then $\phi(S)$ must be a line.  Furthermore, since $S$ maps injectively onto $\phi(S)$, every point of intersection between $\phi(S)$ and the sextic must have even order, so $\phi(C)$ is what is known as a tritangent of $C_F$.  Notice that a section, $S$, and its image under the covering involution, $D-S$, correspond to the same tritangent.

Now consider $NS(\X)$, the N\'eron-Severi lattice of $\X$.  Let $L'$ be the image of $\psi$, the orthogonal projection from $D$.  Let $L\subset L'$ be the sublattice of $NS(\X)$ orthogonal to $D$.  Notice that $L$, unlike $L'$, is integral and even.  Notice also that the index of $L$ in $L'$ is either $1$ or $2$, and must be $2$ if we are to have any elements of norm $-2-\frac{1}{2g-2}=-\frac{5}{2}$, which correspond to sections or tritangents.  The other conditions on the lattices are that we want to minimize the number of vectors $C\in L$ of norm $-2$ and we also have $L\subset NS(\X) \subset H^2(\X,\Z)\cong\II_{3,19}$.  By the following lemma we may as well choose $L$ to have rank $19$.

\begin{lemma}
The maximum number of sections is achieved by a (not necessarily unique) K3 surface with N\'eron-Severi rank $20$.
\end{lemma}

\begin{proof}
It suffices to show that given a K3 surface $\X$ with ample divisor $D$ with $D\cdot D=2$, as long as the rank of $NS(\X)$ is less than $20$ we can increase the rank without decreasing the number of sections or introducing reducible fibers.  Let $c$ be the determinant of $NS(\X)$.  By assumption $T(\X)$ is indefinite, so we can choose a primitive vector $v\in T(\X)$ such that $v\cdot v < -2c^2$.  Consider the lattice $N_1=H^2(\X,\Z)\cap(NS(\X)\otimes\Q + v\Q)$.  We claim that $N_1$ has no vectors of norm $-2$ orthogonal to $D$.

Suppose otherwise.  Let $w=(a/b)v+u\in N_1$ (where $a, b\in\Z$ are relatively prime and $u\in NS(\X)^*$) be such a vector.  Since by assumption there are no such vectors in $NS(\X)$ so we must have $a\ne0$.  Since $v$ and $w$ are both orthogonal to $D$, we find that $u$ must also be orthogonal to $D$ and thus have negative norm.  Since $u\in NS(\X)^*$, $cu\in NS(\X)$; thus, $(da/b)v=cw-cu\in N_1\subset H^2(\X,\Z)$.  Since $v$ was primitive, this implies $b\le c$.  Putting everything together we see that $-2=w\cdot w\le(a^2/b^2)v\cdot v\le(1/c^2)v\cdot v< -2$ which is absurd.

Therefore, by theorem \ref{torelli converse} there exists a K3 surface $\X_1$ such that $NS(\X_1)=N_1$ with $D\in N_1$ ample.  It is clear that the corresponding lattice $L_1'\supset L'$.  Thus, there are at least as many vectors of norm $-\frac{5}{2}$ in $L_1'$ as there are in $L'$.  Now since $\X_1$ has N\'eron-Severi rank larger then $\X$ we are done.
\end{proof}

Furthermore, since gluing a rank $3$ lattice (containing $D$) to $L$ gives us a unimodular lattice, the determinant group must be generated by at most $3$ generators, one of which can be taken to have order $2$ since $D$ has norm $2$.

\subsection{Constructing the Lattice}

We shall construct $L$ (and $NS(\X)$) by starting with a lattice of rank $20$ and taking the sublattice orthogonal to a chosen vector.  As it turns out the following construction works best.

For simplicity we will construct $-L$, the negative of $L$.  Consider the Niemeier lattice with root lattice $A_5^4D_4$.  Then take the part orthogonal to the $D_4$.  This gives us a rank $20$ lattice, $L''$, with root lattice $A_5^4$ and whose determinant group is isomorphic to $\Z/2 \oplus \Z/2$.  

We can also construct $L''$ explicitly.  Let $u\in A_5^*$ be a generator for  $A_5^*/A_5\cong\Z/6$ then $L''$ is generated from $A_5^4$ by the vectors $3u\oplus3u\oplus3u\oplus3u, 0\oplus2u\oplus2u\oplus2u, 2u\oplus2u\oplus4u\oplus0$ \cite[Chap.~16]{cs}.

Since we do not want any vectors of norm $2$, we want a vector that is not orthogonal to any of the roots.  If we identify $A_5=\left\{(x_0,\dots,x_5)\in\Z^6|\sum_i x_i=0\right\}$, then the shortest vectors in $A_5^*$ that are not orthogonal to any root have length $\frac{35}{2}$ and are given by permutations of $v=\frac{1}{2}(-5,-3,-1,1,3,5)$.  Since $v\equiv3u \pmod{A_5}$, we see that $w=v^{\oplus 4}\in L''$.  Note that $N(w)=70$.  Now let $L$ be the negative of the sublattice of  $L''$ orthogonal to $w$.  This $L$ is an even integral lattice with no vectors of norm $-2$ and determinant group $\Z/2 \oplus \Z/2 \oplus \Z/70$.  As it turns out, there is a unique way to glue $L$ to the lattice $\langle D\rangle\cong2\Z$ to get an even integral lattice of signature $(1,19)$ and discriminant $-140$; the transcendental rank 2 lattice one has to glue to get $H^2(\X,\Z)\cong\II_{3,19}$ turns out to have inner-product matrix $\bigl( \begin{smallmatrix}
12&10\\ 10&20
\end{smallmatrix} \bigr)$; this is twice the nonprincipal positive binary quadratic form of discriminant $-35$, which is the only lattice in its genus since the class number of $\Q(\sqrt{-35})$ is $2$.  The corresponding $L'$ has $64$ pairs of elements of norm $-\frac{5}{2}$, giving $64$ tritangents.

I found the construction described above after trying various choices of $L$ and taking the best one.  Although I cannot prove it is the best possible, nevertheless, I feel confident conjecturing it.

\begin{conjecture}
A nonsingular sextic curve in $\P^2$ can have at most $64$ tritangents.  Furthermore, the curve $C_F$ described below is the unique maximum up to the action of $\PGL(3)$.
\end{conjecture}

We next study the symmetry group of the K3 surface, $\X$, corresponding to $L$.  We are interested in the automorphisms that fix the map to $\P^2$, or equivalently, that fix $D$.  Since $D$ is ample and $NS(\X)$ has rank $20$ and thus the transcendental lattice $T(\X)$ has rank $2$, by Theorem \ref{torelli} this automorphism group is the group of isometries of $H^2(\X,\Z)\cong\II_{3,19}$ that fix $D\in H^2(\X,\Z)$ and preserve the projection map onto $T(\X)\otimes\R\cong_\R\C$ up to complex scalar product.  The second condition is equivalent to fixing $T(\X)\subset H^2(\X,\Z)$ and being orientation-preserving on $T(\X)$.  By inspection, the only orientation-preserving automorphism of $T(\X)$ is negation; a corresponding automorphism of $\X$ is the involution $\iota$ corresponding to the double cover $\phi:\X\to\P^2$.  It is easy to see that the action of $\iota$ on $H^2(\X,\Z)$ is given by $E\mapsto (E \cdot D)D - E$, which does indeed act as negation on $T(\X)$ (and $L$).  We now restrict our attention to automorphisms that fix $T(\X)$; this is equivalent to looking at automorphisms of the ramification sextic, $F$, in $\P^2$.  Such an automorphism corresponds to an automorphism of $H^2(\X,\Z)$ that fixes $\langle D\rangle$ and $T(\X)$ element-wise.  Since $H^2(\X,\Z)$ is unimodular, this is equivalent to an automorphism of $L$ that fixes the determinant group $d(L)$ pointwise.  Looking at $L$ we see that this group is the quaternion group $Q_8$.

We next write down this group explicitly.  Let $\AD$ be the anti-diagonal matrix acting on $\Z^6\supset A_5$.   Observe that $\AD$ fixes $A_5$ and sends $v$ to $-v$.  We can now write down the generators
$${\bf i}=\begin{pmatrix} 0&I&0&0 \\ -\AD&0&0&0 \\ 0&0&0&I \\ 0&0&-\AD&0 \end{pmatrix}, \, {\bf j}=\begin{pmatrix} 0&0&0&I \\ 0&0&I&0 \\ 0&-\AD&0&0 \\ -\AD&0&0&0 \end{pmatrix}, {\bf k}=\begin{pmatrix} 0&0&I&0 \\ 0&0&0&-\AD \\ -\AD&0&0&0 \\ 0&I&0&0 \end{pmatrix}.$$
It is easy to write down $6$ vectors in $L$ of norm $-4$ negated by ${\bf i}$, and similarly for ${\bf j}$ and ${\bf k}$.  If $E$ is a divisor in $L$ with self-intersection $-4$, then $D+E$ has self-intersection $-2$ and intersection $2$ with $D$; it therefore $D-E$ corresponds to an effective divisor whose image under $\phi$ is a conic with $6$ tangents $F$.  These conics are nonsingular since the $E$'s are not the sums of $2$ vectors in $L'$ corresponding to tritangent lines.

\subsection{The Action of $Q_8$}

There is a unique faithful action of $Q_8$ on $\P^2(\C)$ coming from the direct sum of the standard action of $Q_8$ on $\C^2$ and the trivial action.  The matrices giving this action are
$${\bf i}=\begin{pmatrix} 0&i&0 \\ i&0&0 \\ 0&0&1 \end{pmatrix}, \, {\bf j}=\begin{pmatrix} 0&-1&0 \\ 1&0&0 \\ 0&0&1 \end{pmatrix}, \, {\bf k}=\begin{pmatrix} i&0&0 \\ 0&-i&0 \\ 0&0&1 \end{pmatrix}.$$
Since the curve $F(X,Y,Z)$ is invariant under this action, the polynomial $F$ must be multiplied by a $1$-dimensional character of $Q_8$; furthermore, that character must be trivial since otherwise $F$ would have no $Z^6$ term and so have a singularity at $[0,0,1]$.  The degree $6$ polynomials invariant under this action are generated by $$X^5Y-Y^5X, \, (X^4+Y^4)Z^2, \, X^2Y^2Z^2, \, {\rm and} \, Z^6.$$  Therefore, $F$ is a linear combination of these functions.  Observe that the coefficients of both $X^5Y-Y^5X$ and $Z^6$ must be nonzero for $F$ to be nonsingular.

Notice also that each of the morphisms ${\bf i}$, ${\bf j}$, ${\bf k}$ has $3$ fixed points in $\P^2(\C)$, corresponding to the eigenvalues of the matrix.  Of these, one is $[0,0,1]$ and the others are on the line $Z=0$ as follows:
\begin{align*}
P_{\bi,1} & =[1,1,0] & P_{\bj,1} & =[i,1,0] & P_{\bk,1} & =[1,0,0] \\
P_{\bi,2} & =[-1,1,0] & P_{\bj,2} & =[-i,1,0] & P_{\bk,2} & =[0,1,0].
\end{align*}
Notice that all the invariant degree $6$ functions, and therefore $F$, vanish on all $6$ of these points; also notice that the tangent lines at these points pass through $[0,0,1]$.  Now consider one of the $6$ conics, tangent to $C_F$ at $6$ points, that is invariant under $\bk$; it is nonsingular, hence the corresponding polynomial has a $Z^2$ term, thus it is invariant under the action of of $\bk$ (not just a scalar multiple of its image under $\bk$).  Thus it must have the form $cXY+Z^2=0$.  In particular the conic passes through $C_F$, hence is tangent to $C_F$, at the points $P_{\bk,m}$.  An analogous argument shows that the $6$ conics invariant under $\bi$ are tangent to $F$ at $P_{\bi,m}$ and similarly for $\bj$.

Now take the pencil of lines in $\P^2$ through the point $[0,0,1]$.  Applying this to the double cover ramified at $F$ gives us a family of genus $2$ curves, $y^2=F(tx,x,1)$, with $65$ pairs of sections (under the hyperelliptic involution) coming from the $64$ tritangents and $1$ corresponding to the base-point, $x=0$.  If we base change our parameterization to the double cover ramified at the lines tangent to $F$ at $P_{\bi,m}$, i.e. the one ramified at $t=\pm1$, we get an additional $12$ pairs of sections, $2$ from each of the $6$ conics that intersect $F$ only at tangent points including $2$ tangent points at $P_{\bi,m}$.  As it turns out there are also $6$ rational quartics in $\P^2$ that again intersect the sextic only at tangent points and furthermore have a node at $[0,0,1]$ whose tangent directions are the lines through $P_{\bi,m}$; each of the quartics gives $2$ additional pairs of sections, giving $24$ extra pairs of sections in all, or $89$ pairs of sections total.  We can do the same construction with $\bj$ to get another cover with $12$ extra pairs of sections.  As it turns out we also get an extra $12$ from rational quartics.  We can then take the fiber product of these two covers to get a $4$-fold cover ramified at $4$ points.  This corresponds to a family of genus $2$ curves parametrized by an elliptic curve with $113$ pairs of sections.

\subsection{Explicit Formulas}

To actually find $F$ explicitly, firstly we find a field over which all of $NS(\X)$ is defined.  Note that $\X$ itself as an abstract K3 surface is defined over $\Q$ since $T(\X)$ is the only lattice in its genus.  Since $\X$ has N\'eron-Severi rank $20$ and $T(\X)$ has discriminant $-140=-2^2\cdot35$ it follows from \cite{in} that $NS(\X)$ is defined over the class-field of the order of height $2$ in $\Q(\sqrt{-35})$.  This class-field is $K=\Q(\sqrt{-35},\sqrt{5},\eta)$ where $\eta$ satisfies $\eta^3 + 3\eta^2 + 5\eta + 1=0$.  Since $71$ is the smallest prime that splits completely in this field we worked over $\F_{71}$.  Since in the equation for the sextic one can scale $X$ and $Y$ independently of $Z$, there are only $2$ free parameters.  Thus it is straightforward to find by exhaustion a sextic of this form in $\P^2(\F_{71})$ that has $64$ tritangents.  Then we used a multivariate Newton iteration to lift the equations for $F$ and its tritangents to the $71$-adics and linear dependence algorithms to interpret the solution over $K$.  Then we used base changes to simplify the coefficients.  The simplest form I got was
\begin{equation*}
\begin{split}
F(X,Y,Z) &= (817022\sqrt{5} - 1925581)X^6 + (1263660\sqrt{5} - 1555968)X^5Y \\ & \quad  + (-1142520\sqrt{5} - 3769035)X^4Y^2 + (3519140\sqrt{5} + 8343710)X^3Y^3 \\ & \quad + (-5279190\sqrt{5} - 11448915)X^2Y^4 + (4070796\sqrt{5} + 9179502)XY^5 \\ & \quad + (-1417660\sqrt{5} - 3164179)Y^6 + \frac{263585\sqrt{5} - 643845}{2}X^4Z^2 \\ & \quad + (214475\sqrt{5} - 322095)X^3YZ^2 + (16305\sqrt{5} - 351765)X^2Y^2Z^2  \\ & \quad + (69845\sqrt{5} + 73095)XY^3Z^2 + (-890\sqrt{5} - 8610)Y^4Z^2 - Z^6.
\end{split}
\end{equation*}

Observe that this $F$ is defined over $\Q(\sqrt{5})$, so we have an action of $G=\Gal(K/\Q(\sqrt{5}))$ on the K3 surface.  We still have our $Q_8$ action; $\bi$ acts as the join of $(\sqrt{-35}/525)M$ and $1$ where $M$ is the $2\times 2$ matrix 
\begin{equation*}
\begin{pmatrix} (27\eta^2 + 54\eta + 45)\sqrt{5} - 66\eta^2 - 51\eta + 82&(-54\eta - 108)\sqrt{5} - 33\eta^2 - 153\eta - 214 \\ (-108\eta^2 - 162\eta - 72)\sqrt{5} + 222\eta^2 + 357\eta + 211&(-27\eta^2 - 54\eta - 45)\sqrt{5} + 66\eta^2 + 51\eta - 82 \end{pmatrix}
\end{equation*}
and $1$; $\bj$ and $\bk$ can be obtained from this $\bi$ by conjugating by elements of $G$.

Under the combined action of $G$ and $Q_8$ there are $3$ orbits of tritangents with representatives
\begin{equation*}
\begin{split}
&((-7\sqrt{-35} + 43)\sqrt{5} + 7\sqrt{-35} - 15)X \\
& + ((\sqrt{-35} - 31)\sqrt{5} + 11\sqrt{-35} - 105)Y + 4\sqrt{5}Z=0,
\end{split}
\end{equation*}
\begin{equation*}
\begin{split}
&((-16\eta^2 - 50\eta - 51)\sqrt{5} - 14\eta^2 - 34\eta - 107)X \\
& + ((-2\eta^2 - 4\eta)\sqrt{5} + 8\eta^2 + 28\eta + 14)Y + 2\sqrt{5}Z=0,
\end{split}
\end{equation*}
and
\begin{equation*}
\begin{split}
&((-26\eta^2 - 81\eta - 130)\sqrt{5} + 76\eta^2 + 201\eta + 308)X \\
& + ((-22\eta^2 - 57\eta - 77)\sqrt{5} + 8\eta^2 + 33\eta + 79)Y + 2\sqrt{5}Z=0.
\end{split}
\end{equation*}
The conics invariant under $\bi$, $\bj$, or $\bk$ have $2$ orbits under the combined action of $G$ and $Q_8$ with representatives
\begin{equation*}
\begin{split}
&((137\eta^2 - 190\eta - 33)\sqrt{5} - 439\eta^2 + 746\eta + 183)X^2 \\
&+ ((164\eta^2 - 478\eta - 162)\sqrt{5} - 124\eta^2 - 334\eta - 162)XY \\
& + ((59\eta^2 + 416\eta + 141)\sqrt{5} + 59\eta^2 + 1064\eta + 357)Y^2 + 4\sqrt{5}Z^2=0
\end{split}
\end{equation*}
and
\begin{equation*}
\begin{split}
&(((2581\sqrt{-35} + 567)\eta^2 + (7014\sqrt{-35} + 462)\eta + 9769\sqrt{-35} - 189)\sqrt{5} \\
& \quad + (-4823\sqrt{-35} + 1155)\eta^2 + (-12138\sqrt{-35} - 6370)\eta - 21259\sqrt{-35} - 105)X^2 \\
& + (((820\sqrt{-35} - 168)\eta^2 + (798\sqrt{-35} + 5502)\eta + 5626\sqrt{-35} + 966)\sqrt{5} \\
& \quad - 4388\sqrt{-35}\eta^2 + (-15498\sqrt{-35} + 12950)\eta - 15094\sqrt{-35} + 630)XY \\
& + (((451\sqrt{-35} - 147)\eta^2 + (3528\sqrt{-35} - 3192)\eta + 1597\sqrt{-35} - 21)\sqrt{5} \\
& \quad + ((103\sqrt{-35} + 105)\eta^2 + (5208\sqrt{-35} - 5320)\eta - 871\sqrt{-35} + 735))Y^2 + 280Z^2=0.
\end{split}
\end{equation*}

We will now construct the above mentioned quartics.  Let $\tilde\X$ be $\X$ blown up at the $2$ points above $[0,0,1]$ thought of as a genus $2$ curve over $\P^1$.  Let $C_1$ be one of the conics fixed by $\bi$.  Let $S_1$ be a double section of $\tilde\X$ above $C_1$ which can be thought of as a point in the Jacobean of the genus $2$ curve; note that $S_1=-\bi^*S_1$.  Let $S_2=\pm\bj^*S_1=\mp\bk^*S_1$ for some choice of $\pm$.  Over a double cover of $\P^1$ these $2$ double sections decompose as the sums of $2$ section $S_1=S_{1,1}+S_{1,2}$, $S_2=S_{2,1}+S_{2,2}$.  Now consider the divisor $S=S_{1,1}-S_{1,2}+S_{2,1}-S_{2,2}$.  Thinking of $S$ as a fiber of the Jacobean, it turns out $S$ can itself be decomposed as $S=Q_1-Q_2$ where $Q_1$ and $Q_2$ are sections over the double cover of $\P^1$.  Then $Q=Q_1+Q_2$ turns out to be a section well-defined over $\P^1$ and its projection to $\P^2$ is one of the quartics.  Since there are $3$ choices of $\{S_1,S_2\}$ and $2$ choices of $\pm$, that gives us $6$ quartics fixed by $\bi$.  The quartics fixed by $\bj$ and $\bk$ are constructed similarly.

\end{document}